\documentclass[12pt,oneside,english]{amsart}

\usepackage[top=2cm,left=2cm,right=2cm,bottom=2cm]{geometry}
\usepackage[english]{babel}
\usepackage{multirow}
\usepackage{graphicx}
\usepackage{amssymb}
\usepackage{yhmath}
\usepackage{verbatim} 
\usepackage{amsmath}
\usepackage{amsthm}
\usepackage{color}
\usepackage{epstopdf} 
\usepackage{wrapfig}
\usepackage[all]{xy}
\usepackage{hyperref}
\usepackage{enumitem} 

\newcommand{\BE}{\begin{equation}}
\newcommand{\EE}{\end{equation}}

\newcommand{\B}{\mathcal{B}}

\newcommand{\Sfe}{\mathbb{S}^3}

\newcommand{\beq}{\begin{eqnarray}}
\newcommand{\eeq}{\end{eqnarray}}

\newcommand{\bev}{\begin{verbatim}}
\newcommand{\eev}{\end{verbatim}}

\newcommand{\vol}{{\rm vol}}

\newtheorem{theorem}{Theorem}

\newcounter{example}[section]

\begin{document}
\title{A topological lower bound for the energy of a unit vector field on a closed hypersurface of the Euclidean space. The 3-dimensional case}
\author{Fabiano G. B. Brito}

\author{Andr\'e O. Gomes}

\author{Adriana V. Nicoli}

\address{Centro de Matem\'atica, Computa\c{c}\~ao e Cogni\c{c}\~ao,
Universidade Federal do ABC,
09.210-170 Santo Andr\'e, Brazil}
\email{fabiano@ime.usp.br}

\address{Dpto. de Matem\'{a}tica, Instituto de Matem\'{a}tica e Estat\'{i}stica,
Universidade de S\={a}o Paulo, R. do Mat\={a}o 1010, S\={a}o Paulo-SP
05508-900, Brazil.}
\email{gomes@ime.usp.br}

\address{Centro de Matem\'atica, Computa\c{c}\~ao e Cogni\c{c}\~ao,
Universidade Federal do ABC,
09.210-170 Santo Andr\'e, Brazil}
\email{dri.nicoli@gmail.com}



\begin{abstract}
In this short note we prove that the degree of the Gauss map $\nu$ of a closed 3-dimensional hypersurface of the Euclidean space is a lower bound for the total bending functional $\mathcal{B}$, introduced by G. Wiegmink. Consequently, the energy functional $E$ introduced by C. M. Wood admits a  topological lower bound.
\end{abstract}
\maketitle

\section{Introduction}
The energy of a unit vector field $\vec{v}$ on a Riemannian compact manifold $M$ is defined by Wood in \cite{Wood}, as the energy of the map $\vec{v}:M\rightarrow T_1M$, where $T_1M$ is the unit tangent bundle equipped with the Sasaki metric, 
\begin{eqnarray}
\label{energy1}
E(\vec{v})=\frac{1}{2}\int_M||\nabla\vec{v}||^2+\frac{m}{2}\vol(M).
\end{eqnarray}
In \cite{Wiegmink}, Wiegmink defines a quantitative measure for the extent to which a unit vector field fails to be parallel with respect to the Levi-Civita connection  $\nabla$ of a Riemannian manifold $M$. This measure is the total bending functional,
\begin{eqnarray}
\label{totalbending}
\mathcal{B}(\vec{v})=\frac{1}{(m-1)\vol(\mathbb{S}^m)}\int_M||\nabla\vec{v}||^2.
\end{eqnarray}
The energy of $\vec{v}$ may be written in terms of the total bending,
\begin{eqnarray}
\label{energy2}
E(\vec{v})=\frac{(m-1)\vol(\mathbb{S}^m)}{2}\mathcal{B}(\vec{v})+\frac{m}{2}\vol(M).
\end{eqnarray}

On the other hand, Brito showed that Hopf flows are absolute minima of the functional $\mathcal{B}$ in $\Sfe$:
\begin{theorem}[Brito, \cite{Brito}]
Hopf vector fields are the unique vector fields on $\Sfe$ to minimize $\mathcal{B}$.
\end{theorem}
Gluck and Ziller proved that Hopf flows are also the unit vector fields of minimum volume, with respect to the following definition of volume,
$$\vol(\vec{v})=\int_M\sqrt{\det(I+(\nabla\vec{v})(\nabla\vec{v})^t)}.$$
\begin{theorem}[Gluck Ziller, \cite{GZ}]
The unit vector fields of minimum volume on $\Sfe$ are precisely the Hopf vector fields, and no others.
\end{theorem}
Reznikov compared this functional to the topology of an Euclidean hypersurface. Let $M$ be a smooth closed oriented immersed hypersurface in $\mathbb{R}^{n+1}$ with the induced metric, and let $\mathcal{S}=\sup_{x\in M}||S_x||=\sup_{x\in M}|\lambda_i(x)|$, where $S_x$ is the second fundamental operator in $T_xM$, and $\lambda_i(x)$ are the principal curvatures.
\begin{theorem}[Reznikov, \cite{Reznikov}]
For any unit vector field $\vec{v}$ on $M$ we have
$$\vol\vec{v}-\vol(M)\geqslant\frac{\vol\mathbb{S}^n}{\mathcal{S}}|\deg(\nu)|,$$
where $\deg(\nu)$ is the degree of the Gauss map $\nu:M\rightarrow\mathbb{S}^n$.
\end{theorem}

Recently, Brito et al, in \cite{Icaro}, discovered a list of curvature integrals for an Euclidean $(n+1)$-dimensional closed hypersurface $M$, 
\begin{equation}
\label{integral-formula}
\int_{M} \eta_k=
\begin{cases}
\deg(\nu){n/2 \choose k/2} {\rm vol}(\mathbb{S}^{n+1}),
  &\mbox{if}\; k \;\mbox{and}\; n \; \mbox{are}\;\mbox{even},\\
0, & \mbox{if}\; k \;\mbox{or}\; n \; \mbox{is}\;\mbox{odd},
\end{cases}
\end{equation}
where the functions $\eta_k$ depend on a unit vector field and on the extrinsic geometry of $M$.

In this note, we take a $3$-dimensional closed Euclidean hipersurface $M$ and relate the total bending (consequently the energy) of a unit vector field $\vec{v}$ to the topology of the manifold $M$, and we prove that 
$${\rm deg}(\nu)\leqslant2\,\widetilde{S}\,\B(\vec{v}),$$
where $\widetilde{S}:=\max_{x\in M}\{||S(\vec{v})_x||\}$, and $S$ is the Weingarten operator of $M$.

\section{Total Bending and Energy of $\vec{v}$}

Let $M$ be a 3-dimensional closed manifold, immersed in the Euclidean space. We may assume that $M$ is oriented, so the normal map $\nu:M\rightarrow\Sfe$,  $\nu(x)=N(x)$, is well defined, where $N$ is a unitary normal field. Let $\nabla$ and $\langle\cdot,\cdot\rangle$ be the Levi-Civita connection and the induced Riemannian metric on $M$, respectively. Let $\vec{v}:M\rightarrow TM$ be a smooth unit vector field on $M$, and take an orthonormal basis $\{e_1,e_2,e_3:=\vec{v}\}$ at each point $x\in M$.   Define the following notation: for $1\leqslant A,B\leqslant3$, $h_{AB}=\langle S(e_A),e_B\rangle$;  for $1\leqslant i,j\leqslant2$, $a_{ij}=\langle\nabla_{e_i}\vec{v},e_j\rangle$, $v_i=\langle\nabla_{\vec{v}}\vec{v},e_i\rangle$.

From the equation \ref{integral-formula}, we know
\begin{equation}
\label{intdeg}
\int_M\eta_2=\deg(\nu)\vol(\Sfe),
\end{equation}
where ${\rm deg}(\nu)$ is the degree of $\nu$ and
\begin{equation}
\label{eta2}
\eta_2=\left|\begin{array}{ccc}
a_{11} & a_{12} & v_1 \\
a_{21} & a_{22} & v_2 \\
h_{31} & h_{32} & h_{33}
\end{array}\right|.
\end{equation}
For more details concerning the functions $\eta_k$, see \cite{Icaro}.
Set $\sigma_1=\left|\begin{array}{cc}
a_{12} & v_1 \\
a_{22} & v_2 \end{array}
\right|, \sigma_2=\left|\begin{array}{cc}
a_{11} & v_1 \\
a_{21} & v_2 \end{array}
\right|$ and $\sigma_3=\left|\begin{array}{cc}
a_{11} & a_{12} \\
a_{21} & a_{22} \end{array}
\right|$. In this case, we have

\begin{eqnarray}
{\rm deg}(\nu)\vol(\Sfe) & = & \int_M\eta_2  =  \int_M\left(h_{31}\sigma_1-h_{32}\sigma_2+h_{33}\sigma_3\right) \nonumber\\
& \leqslant & \left| \int_M\left(h_{31}\sigma_1-h_{32}\sigma_2+h_{33}\sigma_3\right)\right| \nonumber\\
& \leqslant & \int_M\left(|h_{31}||\sigma_1|+|h_{32}||\sigma_2|+|h_{33}||\sigma_3|\right) \nonumber\\
& \leqslant & \frac{1}{2}\int_M\left(|h_{31}|(a_{12}^2+a_{22}^2+v_1^2+v_2^2)+|h_{32}|(a_{11}^2+a_{21}^2+v_1^2+v_2^2)+|h_{33}|(\sum_{i,j=1}^2a_{ij}^2)\right). \nonumber
\end{eqnarray}

Note that $||S(\vec{v})||=\sqrt{h_{31}^2+h_{32}^2+h_{33}^2}$, which implies that $||S(\vec{v})||\geqslant|h_{3A}|$ for $1\leqslant A\leqslant3$. Then

\begin{equation}
\label{quase}
\deg(\nu)\vol(\Sfe)\leqslant \int_M ||S(\vec{v})||\left(\sum_{i,j=1}^2a_{ij}+\sum_{i=1}^2v_i\right).
\end{equation}
By our definition $\widetilde{S}=\max_{x\in M}\{||S(\vec{v})_x||\}$, so
$$\deg(\nu)\vol(\Sfe)\leqslant\widetilde{S}\int_M\left(\sum_{i,j=1}^2a_{ij}+\sum_{i=1}^2v_i\right).$$
On the other hand, by equation \ref{totalbending}, the total bending of $\vec{v}$ may be written as
\begin{equation}
\B(\vec{v})=\frac{1}{2\vol(\Sfe)}\int_M||\nabla\vec{v}||^2.
\end{equation}
Finally 
\begin{equation}
{\rm deg}(\nu)\leqslant2\widetilde{S}\,\B(\vec{v}).
\end{equation}
By equations \ref{energy1} and \ref{energy2}, we deduce the following lower bound for the energy functional:
\begin{eqnarray}
E(\vec{v}) & = & \frac{1}{2}\int_M||\nabla(\vec{v})||^2+\frac{3}{2}\vol(M)\nonumber\\
 & = & \B(\vec{v})\vol(\Sfe)+\frac{3}{2}\vol(M) \nonumber \\
E(\vec{v}) & \geqslant & \frac{1}{2\widetilde{S}}\vol(\Sfe){\rm deg}(\nu)+\dfrac{3}{2}\vol(M) 
\label{lbe}
\end{eqnarray}

Now consider $M=\Sfe$. We are interested in computing the energy of a Hopf vector field $v_H$. By the equation \ref{intdeg}
\begin{equation}
\deg(\nu)\vol(\Sfe)=\int_{\Sfe}\eta_2=\int_{\Sfe}\sigma_3,\nonumber
\end{equation}
since $h_{32}=h_{31}=0$ and $h_{33}=1$. When the canonical immersion of $\Sfe$ in $\mathbb{R}^4$ is considered, we have that $\deg(\nu)=1$. On the other hand, by definition of $E$,
\begin{equation}
\int_{\Sfe}\sigma_3=E(v_H)-\frac{3}{2}\vol(\Sfe),\nonumber
\end{equation}
since $a_{11}=a_{22}=0$ for any Hopf vector field $v_H$. Therefore,
\begin{equation}
E(v_H)=\frac{5}{2}\vol(\Sfe),\nonumber
\end{equation}
as we expected.

\section{Further research}

We intend to extend these preliminary results in two main directions:
\begin{enumerate}
\item Considering the energy of unit vector fields on arbitrary $n$-dimensional hypersurfaces of $\mathbb{R}^{n+1}$.
\item Restricting our attention to totally geodesic flows with integrable normal bundle. This is equivalent of studying the energy of Riemannian foliations of closed hypersurfaces of $\mathbb{R}^3$.
\end{enumerate}

\end{document}